\documentclass[a4paper,12pt]{article}
\usepackage[cp1251]{inputenc}  
\usepackage[english]{babel}
\usepackage{amsmath}
\usepackage{amssymb}
\usepackage{amsfonts}

\def\la  {\lambda}
 \def\b{\beta}    
           \def\ph{\varphi}

\def\be{\begin{equation}}       \def\ba{\begin{array}}
\def\ee{\end{equation}}         \def\ea{\end{array}}

\newtheorem{exi}{Example}

\def\defeq {\stackrel{\mbox{\rm\small def}}{=}}

\newcounter{oldequation}

\begin{document}
\thispagestyle{empty}

\title{Continued fractions and Bessel functions}[Continued fractions and Bessel functions\dots]
    \author{ A.A. Allahverdyan} 
    \newline\hphantom{iii} Adyghe State University,
    \newline\hphantom{iii} Pervomayskaya st., 208, 
    \newline\hphantom{iii} 385000, Maykop, Russia

    \smallskip

    \thanks{\sc Allahverdyan A.A. 
        The Darboux transformations and Bessel functions}
    \thanks{\copyright \ 2019 Allahverdyan A.A.}
    \thanks{\it Recieved ... 2019 }


\maketitle 
\thispagestyle{empty}
\begin{abstract}\small {
Elementary transformations of equations $A\psi=\lambda\psi$ are considered.
The invertibility condition (Theorem 1) is established and similar transformations
of Riccati equations in the case of second order differential operator $A$
are constructed (Theorem 2). Applications to continuous fractions for Bessel functions
 and Chebyshev polynomials are established. It is shown particularly that
 the elementary solutions of Bessel equations are related to a fixed point
 transformations of Riccati equations.
}
\end{abstract}
\medskip

\medskip
\section{Introduction}
Let $A$ be a differential operator $A$ of order $n$
\be\label{formula1}
A=a_0(x)D^{n}_{x}+a_1(x)D^{n-1}_{x} +...+a_n(x), \quad D_x=\frac{d}{dx}.\ee
We consider transformations of this operator defined by substitutions of the form $\hat\psi=(b_0D_x+b_1)\psi$ and their superpositions. In the case $b_0=0$ the transformation is invertible and the operator $A$ in the considered equation transforms in the operator $\hat A$ as follows
 \[\hat A=b_1\circ A\circ b_1^{-1}.\]
The following theorem \cite{cam09} holds true in the general case.
\paragraph{Theorem 1 (on eigenfunctions).} The equation for eigenfunctions $A\psi=\la\psi$, $\la\ne0$ admits an invertible substitution\footnote[1]{Note that the replacement \eqref{formula2} is invertible and its inverse is written by the formula \eqref{formula4}.}
\be\label{formula2}
\hat\psi=(D_x-g)\psi,\quad g=(\log\ph)_x=\frac{\ph_x}{\ph}, \quad A\ph=0.
\ee
First, we prove the following lemma.
\paragraph{Lemma 1.} The differential operator $A$ of order $n>1$ is right divisible by the first order operator $A_1=D_x-g$ iff $g=(\log{\ph})_x$, where $\ph\in\ker A$.\\
$\blacktriangleright$ Let $\phi\in\ker A$ and $g=(\log{\ph})_x$. The formula \eqref{formula1}
\[ A=\sum\limits_{j=0}^{n}{a_j(x)D^{n-j}_{x}}=\tilde A(D_x-f), \quad
\tilde{A}=\sum\limits_{j=0}^{n-1}{\tilde a_j(x)D^{n-(j+1)}_{x}}, \quad g=(\log{\ph})_x,\]
implies that $A(\ph) = 0$ because $(D_x-f)(\ph) = 0$. Then by the substitution $y=\ph\hat y$ we obtain an operator $\hat A$ with the zero coefficient $a_n=0$. Consequently, this polynomial is divisible by $D_x$ iff the initial operator is divisible by $D_x-f$.
$\blacktriangleleft$\\
Now we prove the theorem on eigenfunctions.

$\blacktriangleright$ Note that $A(\ph)=0$ and operator $A$ takes the form \eqref{formula1}:
\[ A=\sum\limits_{j=0}^{n}{a_j(x)D^{n-j}_{x}}.\]
By substituting $A$ in $A\psi=\la\psi$ we find
\be\label{formula3}
\tilde{A}\hat{\psi}=\la\psi\ee
From \eqref{formula3} we have
\be\label{formula4}
\la\psi=a_0\hat{\psi}^{n-1}+a_1\hat{\psi}^{n-2}+...+a_n\hat{\psi}.\ee
Then by the substitution \eqref{formula2} from \eqref{formula4} we obtain
\be\label{formula5}
\la\hat{\psi}=\frac{d}{dx}[a_0\hat{\psi}^{n-1}+a_1\hat{\psi}^{n-2}+...+a_n\hat{\psi}]\ee
 If $\la\ne0$ then the equation \eqref{formula5} and original equation $A\psi=\la\psi$ have the same order, but coefficients in \eqref{formula5} are different.\\
Hence, we have proved that the equation $A\psi=\la\psi$,  $\la\ne0$ admits a substitution $\hat{\psi}=(D_x-g)\psi$ if $g=(\log{\ph})_x$ and this substitution is invertible.\\
$\blacktriangleleft$\\
From this point on, we consider applications of Theorem 1 in the case when $A$ is Euler operator.
\paragraph{Definition 1.} Euler operator has the form
\[A=e^{mt}k(D_t) \]
where $k(D_t)$  is a polynomial in $D_t=\frac{d}{dt}$ with constant coefficients.
\paragraph{Lemma 2.}
If $A=e^{mt}k(D_t)$, $B=e^{nt}z(D_t)$ then a superposition of Euler operators $A$ and $B$ takes the form:
\[A\circ B=e^{(m+n)t}c(D_t),\quad c(D_t)=k(D_t+n)z(D_t).\]
In this case the substitution $\hat\psi=(b_0D_x+b_1)\psi$ becomes an Euler operator of the first order
\be\label{form6}
\hat\psi=e^t(D_t+c)\psi,\quad c\in\mathbb{C}.
\ee
Indeed, $x=e^{-t}$, $dx=-e^{-t}dt$,  therefore $D_x=-e^tD_t$.
\subsection{Second order equations}
Let us consider second order equations and application of Theorem 1 in this case. An operator $A$ can be described as follows
\be\label{7}
A=a_0(x)D^2+a_1(x)D+a_2(x)\ee
Using the substitution $\psi=e^\ph\hat{\psi}$ and assuming that a coefficient of $D$ is equal to zero, we can obtain that a coefficient of $D^2$ is equal to 1, i.e.
\be\label{form8}
A=D^2+q(x).\ee
Then \eqref{form8} takes the form:
\be\label{form9}
A=(D-g)(D+g)\ee
Indeed,
\[A\psi=(D^2+q(x))\psi=\psi''+q(x)\psi.\]
On the other hand,
\[A\psi=(D-g)(D+g)\psi=\psi''+(-g'+g^2)\psi=\psi''+q(x)\psi,\]
where $q(x)+g'+g^2=0$.
\paragraph{Definition 2.} A Riccati equation associated with the equation $A\psi=\la\psi$ is the following equation for the logarithmic derivative
$f=\frac{\psi '}{\psi}$: \be\label{form10} a_0(f'+f^2)+a_1f+a_2=\la.\ee In the particular case that an operator $A$ is of
the form \eqref{form9} the equation \eqref{form10} can be written as follows:
\[f'+f^2+q(x)=\la,\quad f=\frac{\psi'}{\psi}.\]

\section{Bessel equations}
Suppose that an operator $A$ is given by
\be\label{form13}
A=D^2+\dfrac{1}{x}D-\frac{\beta^2}{x^2}.\ee
Then by the substitution  $x=e^{-t}$:
\[D_x=\frac{d}{dx}=-e^{t}\frac{d}{dt}=-e^{t}D_t,\]
one can rewrite the equation $A\psi=\la\psi$ in the following form: \be\label{form14} A\psi=
e^{2t}(D_t^2-\b^2)\psi=e^t(D_t-\b-1)\circ e^t(D_t+\b)\psi=\la\psi.\ee
Note that
\[D_{x}^{2}=(-e^tD_t)\circ (-e^tD_t)=e^{2t}(D_t+D_{t}^{2}).\]
The equation for eigenfunctions of an operator $A$ is $A\psi=\la\psi$. Here an operator $A$ is of the form \eqref{form13}.
The equation considered here is called Bessel equation.\\
By applying Theorem 1 and Lemma 2 to equation \eqref{form14} we obtain \be\label{th1} e^t(D_t+\b)\psi=\hat\psi,\quad
e^t(D_t-\b-1)\hat\psi=\la\psi. \ee As a result, we have that the equation  $A\psi=\la\psi$ takes the form
\be\label{form15} e^{2t}(D_t^2-\hat{\beta}^2)\hat{\psi}=\la\hat{\psi}, \quad \hat\b\defeq \b+1. \ee Rewriting now
equations \eqref{th1} in terms of $f_\b=(\log\psi)_t$ and $ f_{\hat\b }=(\log\hat\psi)_t$ one obtains (see Definition 2):
\be\label{form16f}\begin{gathered} \hat\psi=e^t(f_\b+\b)\psi,\quad \la\psi=e^t( f_{\hat\b }-\hat\b),\\
(f_\b+\b)( f_{\hat\b }-\hat\b)=\la\cdot e^{-2t}=\la\cdot x^2
\end{gathered}\ee
Without loss of generality we put $\la=1$ in the last equation and prove the main theorem.

\paragraph{Theorem 2.} Let  $f=f_\beta$ be a solution of the Riccati equation  $f_t+f^2=\beta^2+x^2$ and the function
$\hat f=f_{\hat{\beta}}$ be defined by the following equation  \be\label{form16} (f+\beta)(\hat{f}-\hat{\beta})=x^2,
\quad \hat\b=\b+1,\ee then this equation states the equivalence of two Riccati equations
\[f_t+f^2=\beta^2+x^2 \Leftrightarrow \hat{f}_t+\hat{f}^2=\hat{\beta}^2+x^2.\]
$\blacktriangleleft$ Let the function $\mu=\mu(t)$ satisfies the differential equation $\mu_t=-2\mu$ and $f_\beta$ be a
solution of the Riccati equation  $f_t+f^2=\beta^2+\mu(t)$. Then the function $f_{\hat{\beta}}$ is defined by the formula
for $\hat{f}$ as follows \be\label{form17} \hat{f}=\frac{\mu}{f+\beta}+\hat{\beta}.\ee By differentiating \eqref{form17}
with respect to $t$: \be\label{form18}
\hat{f}_t=\frac{-2\mu(f+\beta)-f_t\mu}{(f+\beta)^2}=-2\frac{\mu}{f+\beta}-f_t\frac{\mu}{{(f+\beta)^2}}.\ee Note that
\be\label{form19} f_t=\beta^2-f^2+\mu.\ee By substituting \eqref{form19} in \eqref{form18}, \be\label{form20}
\hat{f}_t=-2\frac{\mu}{f+\beta}+\mu\frac{f-\beta}{f+\beta}-\frac{\mu^2}{(f+\beta)^2}.\ee \be\label{form21}
\frac{\mu^2}{(f+\beta)^2}=-\hat{f}_t-2\frac{\mu}{f+\beta}+\mu\frac{f-\beta}{f+\beta}.\ee By squaring both sides of
\eqref{form17} we have
\[\hat{f}^2=\frac{\mu^2}{(f+\beta)^2}+2\hat{\beta}\frac{\mu}{f+\beta}+\hat{\beta}^2.\]
Indeed,
\be\label{form22}
\frac{\mu^2}{(f+\beta)^2}=\hat{f}^2-2\hat{\beta}\frac{\mu}{f+\beta}-\hat{\beta}^2.\ee
Henceforth, from equations \eqref{form21} and \eqref{form22} we can obtain:
\[-\hat{f}_t-2\frac{\mu}{f+\beta}+\mu\frac{f-\beta}{f+\beta}=\hat{f}^2-2\hat{\beta}\frac{\mu}{f+\beta}-\hat{\beta}^2, \]
\[\hat{f}_t+\hat{f}^2=\hat{\beta}^2+\mu\frac{-2+f-\beta+2(\beta+1)}{f+\beta}.\]
So,
\[\hat{f}_t+\hat{f}^2=\hat{\beta}^2+\mu(t).\quad \blacktriangleright \]

\paragraph{Corollary of Theorem 2.} The mapping $A \rightarrow \hat{A}$ defined in Theorem 2 has a fixed point:
\be\label{point} \hat f=f,\quad (\hat\b)^2=\b^2.\ee
 $\blacktriangleleft$ In the case \eqref{point}, \eqref{form16} we find by solving quadratic equations that
\be\label{form23}\begin{gathered} f=f_\pm=\frac{1}{2}\pm x, \quad \b=-\frac12,\quad \hat\b=\frac12,\\
f_t+f^2=\frac14+x^2,\quad (x=e^{-t}).
\end{gathered}\ee
It easy as well to see that $\hat f_\pm=f_\pm\  \blacktriangleright $.
\subsection{Recurrent relations}
 In the case of Chebyshev polynomials $T_n(x)$
 \be\label{chbsh} T_n(x)=\cos n\theta,\quad
x=\cos \theta,\quad D_\theta=-\sqrt{1-x^2}\ D_x\ee we have
 \be\label{chbsh1} T_{n+1}+T_{n-1}=2xT_n,\quad
T_0=1,\quad T_1=x. \ee This recurrent relation totally defines polynomials $T_n(x)$ in $x$-variable and rewriting
\eqref{chbsh1}
 one obtains (Cf. \cite{1866}): \be\label{chbsh2}(f_n-x)(f_{n+1}+x)=-1,\quad f_n=\frac{T_{n-1}}{T_n}+x,\quad f_1(x)=x+\frac1{x}\ee
 whicn looks very similar to \eqref{form16}. Generally speaking,
recurrent relations \eqref{chbsh2} and \eqref{chbsh1} are equivalent. Moreover, by reversing in a certain sense Theorem
2 one may obtain from its proof and \eqref{chbsh2} the second order differential equation
\[D_{x}^{2}T_n=\frac{x\cdot D_x}{1-x^2}T_n-\frac{n^2}{1-x^2}T_n \]
for Chebyshev's polynomials $T_n(x)$. In the case of Bessel functions we can choose, as a basic one, an analog of the linear
recurrent relation \eqref{chbsh1} (see \eqref{form16f} and \cite{Watson}), but the recurrent relation in the ``Riccati''
form \eqref{form16} provides some advantages (Cf \cite{Flajolet}) and yields the formulae \eqref{form23} used below in
order to obtain rational in $x$-variable solutions of the eq. \eqref{form16}.  Introducing a numbering we denote (Cf.
\eqref{form16})
\[ \b_1=\frac12,\quad \beta_{j+1}=\b_j+1,\quad j=1,\ 2,\ \dots \ .\]

\paragraph{Proposition.} Let $\b_j=j-\frac12.$ Then the formula as follows
\be\label{form24} f_{j+1}=\beta_{j+1}+\frac{x^2}{f_j+\beta_j},\quad
 f_1=-x+\frac{1}{2},
\ee provides rational solutions of the Riccati equation of Theorem 2
\[f_2=\frac{3}{2}+\frac{x^2}{1-x}. \]
Similarly determined $f_3, f_4...$ :
\[f_3=\frac{5}{2}+\frac{x^2}{x^2-3x+3};\quad f_4=\frac{7}{2}+\frac{x^2}{6x^2-15x+15}\]
\[\vdots\]

\section*{Conclusion}
Theorem 1 and equation \eqref{form6} reduce the spectral problem $A\psi=\lambda\psi$ with Euler operator $A$ to an algebraic one. This allows us to investigate a generalization of the results of $\S 2$ for higher order Euler operators. Eigenfunctions in this case will provide higher order Bessel functions, but generalization of the continuous fraction approach is not known yet.

\end{document}